# Formulating Oscillator-Inspired Dynamical Systems to Solve Boolean Satisfiability

Mohammad Khairul Bashar, Zongli Lin, Nikhil Shukla*

Department of Electrical & Computer Engineering, University of Virginia, Charlottesville, VA- 22904, USA

*E-mail: ns6pf@virginia.edu

**Abstract**

Dynamical systems can offer a novel non-Boolean approach to computing. Specifically, the natural minimization of energy in the system is a valuable property for minimizing the objective functions of combinatorial optimization problems, many of which are still challenging to solve using conventional digital solvers. In this work, we formulate two oscillator-inspired dynamical systems to solve quintessential computationally intractable problems in Boolean satisfiability (SAT). The system dynamics are engineered such that they facilitate solutions to two different flavors of the SAT problem. We formulate the first dynamical system to compute the solution to the 3-SAT problem, while for the second system, we show that its dynamics map to the solution of the Max-NAE-3-SAT problem. Our work advances understanding of how this physics-inspired approach can be used to address challenging problems in computing.

Dynamical systems offer a unique ‘toolbox’ for solving combinatorial optimization problems [1]-[4]. The intrinsic energy minimization in such systems provides a natural analogue to the minimization of an objective function associated with combinatorial optimization problems [5]. The exploration of new computing paradigms for solving such problems, as in the dynamical system-based approach considered here, is motivated by the fact that computing the solutions to such problems using traditional digital algorithms continues to present a significant challenge [6],[7]. As a case in point, solving Boolean satisfiability (SAT) is an archetypal combinatorial optimization problem that is still considered fundamentally intractable for conventional digital hardware. The SAT problem is defined as the challenge of evaluating if there exists a Boolean assignment for the variables in a given Boolean expression (in the conjunctive normal form) that would make the expression TRUE. Besides being the first known NP-complete problem [8], the SAT problem is considered particularly relevant since many practical combinatorial optimization problems can be easily reduced to the solution of the SAT problem. Here, we specifically consider the case of the 3-SAT problem, a constrained but NP-complete version of SAT, where each clause contains no more than three literals.

In this work, we formulate and analyze two oscillator-inspired dynamical systems and show that their dynamics can be directly used to compute solutions to the 3-SAT (System I) and the Max-NAE-3-SAT problem (System II). The Not-all-Equal (NAE)-SAT problem is an NP-complete variant of the SAT problem which imposes the additional constraint that every clause must contain a literal that is true and another literal that is false; the Max-NAE-SAT problem is the optimization version of the problem where the objective is to maximize the number of clauses that meet this constraint. We note that Ercsey-Ravasz *et al.* [9] proposed an analog computational model for solving the SAT problem which was formulated using non-oscillating (analog) variables; further, in our prior work, we have also proposed computational models for many combinatorial problems (e.g., NAE-SAT, integer factorization among others) with non-oscillating analog variable [10]. While we draw many important insights from these works, our effort here is fundamentally different in that our dynamical systems use oscillating (analog) variables, and consequently, exhibits a different set of dynamics. We would also like to point out that there have been multiple prior works that have explored the formulation of oscillator-based computational models for solving combinatorial optimization problems such as Maximum Cut [11]-[15], Maximum Independent Set [16],[17], Graph coloring [18]-[21] and Max-K-Cut [22] among others. However, all these problems, unlike Boolean satisfiability, have objective functions with quadratic degree [23]. Subsequently, the oscillator-based computational models can be developed using the Kuramoto framework which cannot be *directly* applied here.

**System I:**

To formulate System I, we represent every variable $x_i$ in the Boolean expression using an analog variable $\alpha_i$, where $x_i = \frac{1+\cos(t+\alpha_i)}{2}$, which can be considered as a level shifted oscillator; the oscillator angular frequency ($\omega$) is assumed to be $\omega = 1$ in this theoretical analysis. The relationship between $x_i$ and $\alpha_i$ $\left(x_i = \frac{1+\cos(t+\alpha_i)}{2}\right)$ is defined such that the maximum (or minimum) value of the analog variable equals the Boolean assignment for $x_i \in \{0,1\}$, respectively. For each clause $C_m$, we define $K_{m,osc}(t,\alpha) = \prod_{i=1}^{N}\left(1 - \left(\frac{1+c_{mi}\cos(t+\alpha_i)}{2}\right)\right)$, where $c_{mi} = 1(-1)$, if the $i^{th}$ variable appears in the $m^{th}$ clause in the normal (negated) form; $c_{mi} = 0$, if the variable is absent from the $m^{th}$ clause; $\alpha = [\alpha_1\ \alpha_2 \ldots\ldots \alpha_N]$; $N$ is the number of variables in the SAT problem. It can be observed that $K_{m,osc}(t, \alpha) = 0$, if and only if the clause is satisfied. We define the dynamical system: $(-\nabla_\alpha V)_i = 1 + \frac{d\alpha_i}{dt}$. The energy function for the system is defined as:

$$V = \sum_{m=1}^{M} A\left(K_{m,osc}(t, \alpha)\right)^2 \tag{1}$$

Here, $M$ is the total number of clauses in the problem. $V = 0$ when all the clauses are satisfied, and consequently, corresponds to the solution of the SAT problem (if the problem is satisfiable). To evaluate the temporal evolution of the system energy, we calculate $\frac{dV}{dt}$, which is given by:

$$\frac{dV}{dt} = \sum_{i=1}^{N} \left(\frac{\partial V}{\partial \alpha_i}\right)\left(\frac{d\alpha_i}{dt}\right) + \frac{\partial V}{\partial t} \tag{2}$$

Using equation (1) and the definition of $K_{m,osc}(t, \alpha)$, we can calculate $\frac{\partial V}{\partial t}$ as,

$$\begin{aligned} \frac{\partial V}{\partial t} &= \sum_{m=1}^{M} \left( 2AK_{m,osc}(t, \alpha)\ \frac{\partial\left(K_{m,osc}(t, \alpha)\right)}{\partial t} \right) \\ &= \sum_{m=1}^{M} \left( 2AK_{m,osc}(t, \alpha).\ \left[ \sum_{i=1}^{N} \frac{c_{mi} K_{m,osc}(t, \alpha)}{1 - c_{mi}\cos(t + \alpha_i)} \sin(t + \alpha_i) \right] \right) \\ &= \sum_{i=1}^{N} \left[ \sum_{m=1}^{M} \left( 2AK_{m,osc}(t, \alpha)\ \frac{c_{mi} K_{m,osc}(t, \alpha)}{1 - c_{mi}\cos(t + \alpha_i)} \sin(t + \alpha_i) \right) \right] \end{aligned} \tag{3}$$

Further, $\frac{\partial V}{\partial \alpha_i}$ can be calculated as,

$$\frac{\partial V}{\partial \alpha_i} = \sum_{m=1}^{M} \left( 2AK_{m,osc}(t,\alpha)\ \frac{\partial \left( K_{m,osc}(t,\alpha) \right)}{\partial \alpha_i} \right)$$

$$= \sum_{m=1}^{M} \left( 2AK_{m,osc}(t,\alpha)\ \frac{c_{mi} K_{m,osc}(t,\alpha)}{1 - c_{mi}\cos(t+\alpha_i)} \sin(t+\alpha_i) \right) \tag{4}$$

Substituting equation (4) into equation (3), $\frac{\partial V}{\partial t}$ can be expressed as,

$$\frac{\partial V}{\partial t} = \sum_{i=1}^{N} \frac{\partial V}{\partial \alpha_i} \tag{5}$$

By substituting the expression for $\frac{\partial V}{\partial t}$ from equation (5) into equation (2), $\frac{dV}{dt}$ can be calculated as,

$$\frac{dV}{dt} = \sum_{i=1}^{N} \left(\frac{\partial V}{\partial \alpha_i}\right)\left(\frac{d\alpha_i}{dt}\right) + \frac{\partial V}{\partial t} = \sum_{i=1}^{N} \left(\frac{\partial V}{\partial \alpha_i}\right)\left(\frac{d\alpha_i}{dt}\right) + \sum_{i=1}^{N} \frac{\partial V}{\partial \alpha_i} = \sum_{i=1}^{N} \left(\frac{\partial V}{\partial \alpha_i}\right)\left(1 + \frac{d\alpha_i}{dt}\right) \tag{6a}$$

Further, utilizing the system dynamics $(-\nabla_\alpha V)_i = 1 + \frac{d\alpha_i}{dt}$ (defined above), equation (6a) can be expressed as,

$$\frac{dV}{dt} = \sum_{i=1}^{N} \left(\frac{\partial V}{\partial \alpha_i}\right)\left(1 + \frac{d\alpha_i}{dt}\right) = -\sum_{i=1}^{N} \left(1 + \frac{d\alpha_i}{dt}\right)\left(1 + \frac{d\alpha_i}{dt}\right) = -\sum_{i=1}^{N} \left(1 + \frac{d\alpha_i}{dt}\right)^2 \tag{6b}$$

It can be observed from equation 6(b) that $V$ is a decreasing function with time since $\frac{dV}{dt} \leq 0$. Consequently, this implies that the corresponding system dynamics will evolve to reduce system energy $(V)$.

In order to formulate the system dynamics $\frac{d\alpha_i}{dt}$, we express $\frac{dV}{dt}$ as,

$$\frac{dV}{dt} = \sum_{m=1}^{M} \left( 2AK_{m,osc}(t,\alpha) \frac{d\left(K_{m,osc}(t,\alpha)\right)}{dt} \right) \tag{7a}$$

$$\frac{dV}{dt} = \sum_{m=1}^{M} \left( 2AK_{m,osc} \sum_{i=1}^{N} \left( c_{mi} \frac{K_{m,osc}}{1 - c_{mi}\cos(t+\alpha_i)} \sin(t+\alpha_i). \left(1 + \frac{d\alpha_i}{dt}\right) \right) \right) \tag{7b}$$

$$\frac{dV}{dt} = \sum_{i=1}^{N} \left( \sum_{m=1}^{M} \left( 2AK_{m,osc} c_{mi} \frac{K_{m,osc}}{1 - c_{mi}\cos(t+\alpha_i)} \right) \right) \sin(t+\alpha_i). \left(1 + \frac{d\alpha_i}{dt}\right) \tag{7c}$$

Equating (6b) and (7c), we get

$$-\left(1 + \frac{d\alpha_i}{dt}\right) = \sum_{m=1}^{M} \left( 2AK_{m,osc}(t,\alpha) \left[ \frac{c_{mi} K_{m,osc}(t,\alpha)}{1 - c_{mi}\cos(t+\alpha_i)} \right] \sin(t+\alpha_i) \right) \tag{8a}$$

Equation (8a) can be rewritten as,

$$\frac{d\alpha_i}{dt} = -\left( \sum_{m=1}^{M} \left( 2AK_{m,osc}(t,\alpha) \left[ \frac{c_{mi} K_{m,osc}(t,\alpha)}{1 - c_{mi}\cos(t+\alpha_i)} \right] \sin(t+\alpha_i) \right) + 1 \right) \tag{8b}$$

The right-hand side in equation (8b) is 2π periodic in time. At steady state, $V = 0$; $\frac{d\alpha_i}{dt} = -1$, and consequently, $\alpha_i = -t + c_i$, where $c_i$ is a constant. The system is designed such that the out-of-phase feedback essentially 'cancels' out the oscillations when the system achieves the ground state energy. This corresponds to all the clauses being satisfied (if the problem is satisfiable). Figure 1 illustrates the system dynamics for a representative SAT problem. Details of the simulation framework have been discussed in Appendix I.

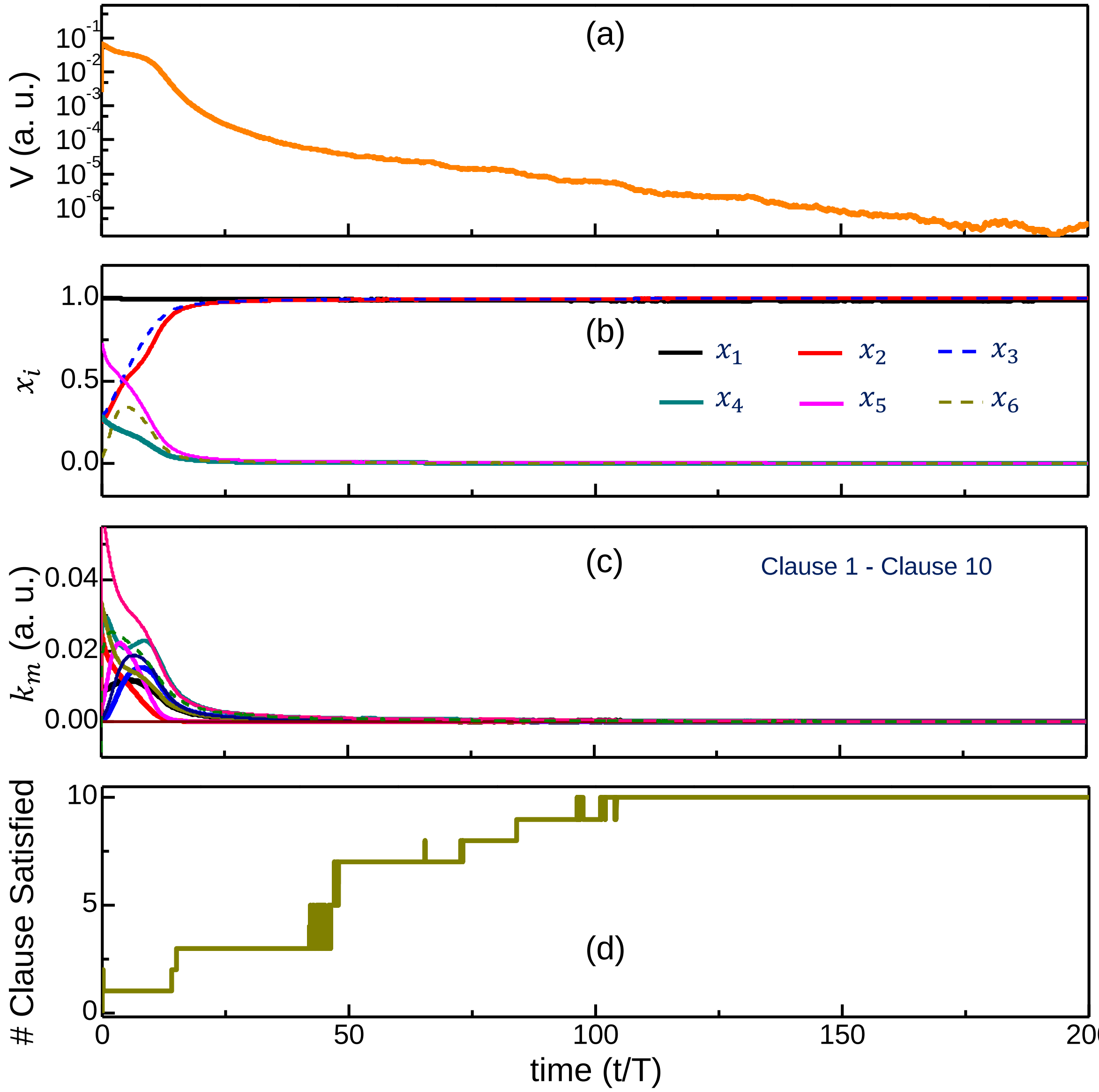

**Fig. 1.** Evolution of (a) $V$, (b) $x_i$, (c) $k_m$, and (d) No. of clauses satisfied with time for an illustrative 3-SAT problem with 6-variables and 10-clauses that is computed using the System I dynamics. Here, $\omega = 2\pi$ is used such that $T = 1$. The simulation is performed using a stochastic differential equation framework (details in Appendix I).

**System II:**

For this implementation, we formulate the system dynamics as: $(-\nabla_\alpha E)_i = \frac{d\alpha_i}{dt}$, where $E$ is the potential energy function of the system. In contrast to the prior approach, here, we will first define the system dynamics, and subsequently, aim to show that there exists a Lyapunov (energy) function which can directly be mapped to the solution to the Max-NAE-3-SAT problem. We consider a system whose dynamics are defined by:

$$\frac{d\alpha_i}{dt} = \sin(t+\alpha_i)\left(-\sum_{m=1}^{M}\left(2AK_{m,osc}(t,\alpha)\left[\frac{c_{mi}K_{m,osc}(t,\alpha)}{1-c_{mi}\cos(t+\alpha_i)}\right]\right)\right) \tag{9}$$

$$-\sin(2t+2\alpha_i).\left[A_s\cos(2t)\right]$$

$$\equiv \chi\big(t+\alpha_i(t)\big).B_i(t) + \chi\big(2t+2\alpha_i(t)\big).B^{(2)}(t)$$

Equation (9) can be interpreted as a (sinusoidal) oscillator under perturbation ($B_i(t)$), and second harmonic signal injection $B^{(2)}(t) \equiv A_s\cos(2t)$ which helps binarize the phases to (0, π) [24],[25], as illustrated further on. $\chi(t+\alpha_i)$ and $\chi(2t+2\alpha_i)$ are the first and the second harmonics of the perturbation projection vectors (PPVs) of the oscillator, respectively. $A$ and $A_s$ are positive constants. It can be observed that the dynamics described in equation (9) are a modified version of the dynamics derived in equation (8b) for System I, and essentially help us formulate the dynamics of System II. However, it must be emphasized here that we do not use the potential energy function $V$ defined for System I since it does not decrease monotonically for the System II dynamics. Instead, using the dynamics described above, we will formulate a new energy function $E$ whose ground state maps to the solution to the Max-NAE-3-SAT problem.

To define $E$, we first reformulate equation (9) in terms of the relative phase difference. Substituting the definition of $K_{m,osc}(t,\alpha)$, equation (9) can be rewritten as,

$$\frac{d\alpha_i}{dt} = -A\sin(t+\alpha_i).\sum_{m=1}^{M}\left[c_{mi}\left(\prod_{j=1;j\neq i}^{N}\left(\frac{1-c_{mj}\cos(t+\alpha_j)}{2}\right)\right)^2\left(\frac{1-c_{mi}\cos(t+\alpha_i)}{2}\right)\right] \qquad (10)$$
$$-\sin(2t+2\alpha_i).[A_s\cos(2t)]$$

Expanding equation (10),

$$\frac{d\alpha_i}{dt} = -\frac{A}{2}\left[\sum_{m=1}^{M}\left(c_{mi}\sin(t+\alpha_i)\left(\prod_{j=1;j\neq i}^{N}\left(\frac{1-c_{mj}\cos(t+\alpha_j)}{2}\right)\right)^2\right)\right.$$
$$\left.-\sum_{m=1}^{M}\left(\frac{1}{2}c_{mi}^2\sin(2(t+\alpha_i))\left(\prod_{j=1;j\neq i}^{N}\left(\frac{1-c_{mj}\cos(t+\alpha_j)}{2}\right)\right)^2\right)\right] \qquad (11)$$
$$-\sin(2t+2\alpha_i).[A_s\cos(2t)]$$

Further, using trigonometric identities to express all the product terms in $\left(\prod_{j=1;j\neq i}^{N}\left(\frac{1-c_{mj}\cos(t+\alpha_j)}{2}\right)\right)^2$ as the sum of $\cos(.)$ terms, we rewrite the expression as,

$$\sum_{\mu_N=-2}^{2}\cdots\sum_{\mu_2=-2}^{2}\sum_{\mu_1=-2}^{2}C_{\mu_1,\mu_2\ldots\mu_N;\neq\mu_i}\cos\left(\left(\sum_{j=1;\, j\neq i}^{N}|c_{mj}|\mu_j\right)t+\sum_{j=1;\, j\neq i}^{N}|c_{mj}|\mu_j\alpha_j\right)$$

Using the approach described by Wang *et al.* [26], a differential equation such as equation (11) can be formulated as a Multi-time Partial Differential Equation (MPDE), wherein the fundamental oscillation is assumed to happen in fast time $t_1$ while the phases evolve in slow time $t_2$. Subsequently, equation (11) can then be approximated as,

$$\frac{d\alpha_i}{dt} = -A\sum_{m=1}^{M}\left(\sum_{\mu_N=-2}^{2}\cdots\sum_{\mu_2=-2}^{2}\sum_{\mu_1=-2}^{2} c_{mi}\ Q_1 C^1_{\mu_1,\mu_2\ldots\mu_N;\neq\mu_i}\sin\left(\alpha_i - \sum_{j=1;\ j\neq i}^{N}\left|c_{mj}\right|\mu_j\alpha_j(t)\Bigg|_{Q_1}\right)\right)$$

$$+$$

$$A\sum_{m=1}^{M}\left(\sum_{\mu_N=-2}^{2}\cdots\sum_{\mu_2=-2}^{2}\sum_{\mu_1=-2}^{2} c^2_{mi}\ Q_2 C^2_{\mu_1,\mu_2\ldots\mu_N;\neq\mu_i}\sin\left(2\alpha_i - \sum_{j=1;\ j\neq i}^{N}\left|c_{mj}\right|\mu_j\alpha_j(t)\Bigg|_{Q_2}\right)\right) \quad (12)$$

$$-A_{s1}\sin(2\alpha_i)$$

Here, $Q_1 = 1$ when $\sum_{j=1;\ j\neq i}^{N}\left|c_{mj}\right|\mu_j = 1$ else $Q_1 = 0$; $Q_2 = 1$ when $\sum_{j=1;\ j\neq i}^{N}\left|c_{mj}\right|\mu_j = 2$ else $Q_2 = 0$. Additional details regarding the derivation of equation (12) can be found in Appendix II. Remarkably, there is a Lyapunov function $E\big(\alpha(t)\big)$ which can be defined for these dynamics as,

$$E\big(\alpha(t)\big) = \sum_{i=1}^{N}\left[-A\sum_{m=1}^{M}\sum_{\mu_N=-2}^{2}\cdots\sum_{\mu_2=-2}^{2}\sum_{\mu_1=-2}^{2} c_{mi}\ Q_1\ C^1_{\mu_1,\mu_2\ldots\mu_N;\neq\mu_i}\cos\left(\alpha_i(t) - \sum_{j=1;\ j\neq i}^{N}\left|c_{mj}\right|\mu_j\alpha_j(t)\Bigg|_{\sum_{j=1;\ j\neq i}^{N}\left|c_{mj}\right|\mu_j=1}\right)\right]$$

$$+\sum_{i=1}^{N}\left[\frac{A}{2}\sum_{m=1}^{M}\sum_{\mu_N=-2}^{2}\cdots\sum_{\mu_2=-2}^{2}\sum_{\mu_1=-2}^{2} c^2_{mi}\ Q_2\ C^2_{\mu_1,\mu_2\ldots\mu_N;\neq\mu_i}\cos\left(2\alpha_i(t) - \sum_{j=1;\ j\neq i}^{N}\left|c_{mj}\right|\mu_j\alpha_j(t)\Bigg|_{\sum_{j=1;\ j\neq i}^{N}\left|c_{mj}\right|\mu_j=2}\right)\right] - \frac{A_{s1}}{2}\cos\big(2\alpha_i(t)\big) \quad (13)$$

Unlike $V$ (defined for System I), $E\big(\alpha(t)\big)$ is defined in terms of relative phase difference (and not in terms of the absolute phase). To show that $E\big(\alpha(t)\big)$ is a decreasing function in time i.e., $\frac{dE(\alpha(t))}{dt} \leq 0$, we express $\frac{dE(\phi(t))}{dt} = \frac{dE(\alpha(t))}{d\alpha_i(t)} . \frac{d\alpha_i(t)}{dt}$, where $\frac{dE(\alpha(t))}{d\alpha_i(t)}$ can be calculated as,

$$\frac{\partial E(\alpha(t))}{\partial \alpha_i(t)} =$$

$$A \sum_{m=1}^{M} \sum_{\mu_N=-2}^{2} .. \sum_{\mu_2=-2}^{2} \sum_{\mu_1=-2}^{2} c_{mi}\ Q_1\ C^1_{\mu_1,\mu_2 \dots \mu_N; \neq \mu_i} \sin\left( \alpha_i - \sum_{j=1;\ j\neq i}^{N} |c_{mj}| \mu_j \alpha_j(t) \Bigg|_{Q_1} \right)$$

$$+$$

$$-\frac{2A}{2} \sum_{m=1}^{M} \sum_{\mu_N=-2}^{2} .. \sum_{\mu_2=-2}^{2} \sum_{\mu_1=-2}^{2} c^2_{mi}\ Q_2\ C^2_{\mu_1,\mu_2 \dots \mu_N; \neq \mu_i} \sin\left( 2\alpha_i - \sum_{j=1;\ j\neq i}^{N} |c_{mj}| \mu_j \alpha_j(t) \Bigg|_{Q_2} \right) + A_{s1} \sin(2\alpha_i) \tag{14}$$

$$\equiv -\frac{d\alpha_i(t)}{dt}$$

Thus,

$$\frac{\partial E(\alpha(t))}{\partial \alpha_i(t)} = -\frac{d\alpha_i(t)}{dt} \tag{15}$$

It can be observed that equation (15) represents the system dynamics described earlier. Subsequently,

$$\frac{dE(\alpha(t))}{dt} = \sum_{i=1}^{N} \left[ \left( \frac{\partial E(\alpha(t))}{\partial \alpha_i(t)} \right) . \left( \frac{d\alpha_i(t)}{dt} \right) \right] \tag{16}$$

$$= -\sum_{i=1}^{N} \left[ \left( \frac{d\alpha_i(t)}{dt} \right)^2 \right] \leq 0 \tag{17}$$

Equation (17) reveals that $E(\alpha(t))$ is decreasing in time.

While equation (13) represents a general form, we will specifically define the energy $E$ for the case when each clause contains exactly 3 literals, and subsequently, show that its ground state can be used to find the solution of the NAE-3-SAT problem. When a clause contains 3 literals (corresponding to variables $i, j, k$), $E$ can be expressed as,

$$E(\alpha)=\sum_{i=1}^{N}\Bigg[\pi A.2^{-2N+1}\sum_{\substack{m=1;i\neq j\neq k;c_{mi}\neq 0\\ c_{mj}\neq 0,c_{mk}\neq 0}}^{M}\Bigg(2c_{mi}c_{mj}\left(1+\frac{1}{2}c_{mk}^{2}\right)\cos(\alpha_i-\alpha_j)$$
$$+2c_{mi}c_{mk}\left(1+\frac{1}{2}c_{mj}^{2}\right)\cos(\alpha_i-\alpha_k)+\frac{1}{2}c_{mi}c_{mj}c_{mk}^{2}\cos(\alpha_i+\alpha_j-2\alpha_k)$$
$$+\frac{1}{2}c_{mi}c_{mk}c_{mj}^{2}\cos(\alpha_i+\alpha_k-2\alpha_j)+\frac{1}{8}c_{mi}^{2}c_{mk}^{2}\left(1+\frac{1}{2}c_{mj}^{2}\right)\cos(2\alpha_i-2\alpha_k) \quad (18)$$
$$+\frac{1}{2}c_{mi}^{2}c_{mj}c_{mk}\cos(2\alpha_i-\alpha_j-\alpha_k)$$
$$+\frac{1}{8}c_{mi}^{2}c_{mj}^{2}\left(1+\frac{1}{2}c_{mk}^{2}\right)\cos(2\alpha_i-2\alpha_j)\Bigg)\Bigg]-\sum_{i=1}^{N}\frac{\pi A_s}{2}\cos(2\alpha_i)$$

The details of this derivation are shown in Appendix III. We note that if a clause contains literals corresponding to only one or two distinct variables, $i \neq j \neq k$ constraint will not be imposed for that specific clause in equation (18). The specific nature of the arguments of the $\cos(.)$ functions shown in equation (18) arise from the characteristics of the cross-correlation operation performed in equation (12). The corresponding dynamics associated with equation (18) can be defined as,

$$\frac{d\alpha_i}{dt} = \pi A. 2^{-2N+1} \sum_{\substack{m=1; i\neq j\neq k; c_{mi}\neq 0 \\ c_{mj}\neq 0, c_{mk}\neq 0}}^{M} \Big[ 2c_{mi}c_{mj}\left(1+\frac{1}{2}c_{mk}^2\right)\sin(\alpha_i-\alpha_j)$$
$$+ 2c_{mi}c_{mk}\left(1+\frac{1}{2}c_{mj}^2\right)\sin(\alpha_i-\alpha_k) + \frac{1}{2}c_{mi}c_{mj}c_{mk}^2\sin(\alpha_i+\alpha_j-2\alpha_k)$$
$$+\frac{1}{2}c_{mi}c_{mk}c_{mj}^2\sin(\alpha_i+\alpha_k-2\alpha_j) + \frac{1}{4}c_{mi}^2c_{mk}^2\left(1+\frac{1}{2}c_{mj}^2\right)\sin(2\alpha_i-2\alpha_k) \tag{19}$$
$$+\frac{1}{4}c_{mi}^2c_{mj}^2\left(1+\frac{1}{2}c_{mk}^2\right)\sin(2\alpha_i-2\alpha_j)$$
$$+ c_{mi}^2c_{mj}c_{mk}\sin(2\alpha_i-\alpha_j-\alpha_k)\Big] - \pi A_s\sin(2\alpha_i)$$

The second harmonic injection signal $-\sum_{i=1}^{N}\frac{\pi A_s}{2}\cos(2\alpha_i)$ (for an appropriate injection strength $A_s$) essentially lowers the energy of the system corresponding to $\alpha \in \{0, \pi\}$, since the minimization of $-\sum_{i=1}^{N}\frac{\pi A_s}{2}\cos(2\alpha_i)$ to $-\sum_{i=1}^{N}\frac{\pi A_s}{2}$ forces the oscillators to take these binary phase values; this concept was also exploited in designing oscillator-based Ising machines [26]. Thus, when the system achieves ground state, each $2\alpha$ term in equation (18) induces a phase difference of $0$ or $2\pi$, and hence, the arguments of the corresponding $\cos(\alpha_i+\alpha_j-2\alpha_k)$ terms can be simplified to $\cos(\alpha_i+\alpha_j)$. Further, $\cos(2\alpha_i-2\alpha_j)$ will take constant values at these specific phase points (represented as $C$). Additionally, $c_{mi}^2 = c_{mj}^2 = c_{mk}^2 = 1$. Thus, at these discrete phase points, $E(\alpha)$ for a problem in which each clause consists of three literals can be reduced to,

$$E(\alpha) = \pi A. 2^{-2N+1} \sum_{i=1}^{N} \left[ \sum_{\substack{m=1; i \neq j \neq k; c_{mi} \neq 0 \\ c_{mj} \neq 0, c_{mk} \neq 0}}^{M} \left( 3c_{mi}c_{mj}\cos(\alpha_i - \alpha_j) + 3c_{mi}c_{mk}\cos(\alpha_i - \alpha_k) \right. \right.$$

$$+ \frac{1}{2}c_{mi}c_{mj}\cos(\alpha_i + \alpha_j) + \frac{1}{2}c_{mi}c_{mk}\cos(\alpha_i + \alpha_k) \tag{20}$$

$$\left. \left. + \frac{1}{2}c_{mj}c_{mk}\cos(\alpha_j + \alpha_k) \right) \right] + C - \sum_{i=1}^{N} \frac{\pi A_s}{2}\cos(2\alpha_i)$$

Rearranging equation (20),

$$E(\alpha) = \pi A. 2^{-2N+1} \sum_{\substack{m=1; i \neq j \neq k; c_{mi} \neq 0 \\ c_{mj} \neq 0, c_{mk} \neq 0}}^{M} \left[ \sum_{i=1}^{N} \left( 3c_{mi}c_{mj}\cos(\alpha_i - \alpha_j) + 3c_{mi}c_{mk}\cos(\alpha_i - \alpha_k) \right. \right.$$

$$+ \frac{1}{2}c_{mi}c_{mj}\cos(\alpha_i + \alpha_j) + \frac{1}{2}c_{mi}c_{mk}\cos(\alpha_i + \alpha_k)$$

$$\left. \left. + \frac{1}{2}c_{mj}c_{mk}\cos(\alpha_j + \alpha_k) \right) \right] + C - \sum_{i=1}^{N} \frac{\pi A_s}{2}\cos(2\alpha_i)$$

$$= \sum_{m=1}^{M} \beta_m(\alpha_i, \alpha_j, \alpha_k) + C - \sum_{i=1}^{N} \frac{\pi A_s}{2}\cos(2\alpha_i) = \sum_{m=1}^{M} \beta_m(\alpha_i, \alpha_j, \alpha_k) + C - C_s \tag{21}$$

Both $C$ and $C_s \left(= \sum_{i=1}^{N} \frac{\pi A_s}{2}\cos(2\alpha_i)\right)$ are constants at the phase points, $\alpha \in \{0, \pi\}$. Consequently, $E(\alpha)$ is minimized when $\sum_{m=1}^{M} \beta_m(\alpha_i, \alpha_j, \alpha_k)$ is minimum. Now, for a single clause consisting of 3 literals corresponding to 3 variables $x_i$, $x_j$, $x_k$ (here, $x_i$, $x_j$, $x_k$ can appear in normal or negated form in the clause), $\beta_m(\alpha_i, \alpha_j, \alpha_k)$ can be written as,

$$\beta_m(\alpha_i, \alpha_j, \alpha_k) = \pi A. 2^{-2N+1} \Big[ c_{mi} c_{mj} \left( 6\cos(\alpha_i - \alpha_j) + \frac{3}{2}\cos(\alpha_i + \alpha_j) \right) + c_{mj} c_{mk} \left( 6\cos(\alpha_j - \alpha_k) + \frac{3}{2}\cos(\alpha_j + \alpha_k) \right) + c_{mk} c_{mi} \left( 6\cos(\alpha_k - \alpha_i) + \frac{3}{2}\cos(\alpha_k + \alpha_i) \right) \Big] + C - C_s \quad (22a)$$

$$\beta_m(\alpha_i, \alpha_j, \alpha_k) = \pi A. 2^{-2N+1} [T_{ij} + T_{jk} + T_{ki}] + C - C_s \quad (22b)$$

where

$$T_{ij} = c_{mi} c_{mj} \left( 6\cos(\alpha_i - \alpha_j) + \frac{3}{2}\cos(\alpha_i + \alpha_j) \right) \quad (23)$$

Equation (22b) reveals that $\beta_m(\alpha_i, \alpha_j, \alpha_k)$ is minimum when $[T_{ij} + T_{jk} + T_{ki}]$ is minimum.

At the phase points $\alpha_i, \alpha_j, \alpha_k \in \{0, \pi\}$, $T_{ij}, T_{jk}, T_{ki}$ and $[T_{ij} + T_{jk} + T_{ki}]$ are binary in nature and exhibits the property that $[T_{ij} + T_{jk} + T_{ki}]$, and thus $\beta_m(\alpha_i, \alpha_j, \alpha_k)$, is minimized when $(x_i \oplus x_j) \vee (x_j \oplus x_k) \vee (x_k \oplus x_i) = 1$. This is illustrated in the following paragraph. However, first, we simplify $(x_i \oplus x_j) \vee (x_j \oplus x_k) \vee (x_k \oplus x_i)$ as,

$$(x_i \oplus x_j) \vee (x_j \oplus x_k) \vee (x_k \oplus x_i) = (x_i \overline{x_j} \vee \overline{x_i} x_j) \vee (x_j \overline{x_k} \vee \overline{x_j} x_k) \vee (x_k \overline{x_i} \vee \overline{x_k} x_i) = (x_i \vee x_j \vee x_k). (\overline{x_i} \vee \overline{x_j} \vee \overline{x_k}) \quad (24)$$

Remarkably, equation (24) corresponds to a clause of the NAE-3-SAT problem. The terms within the first parentheses in equation (24) implement the standard SAT constraint while the terms in the second parentheses implement the constraint that at least one literal must be false. Here, we again emphasize that $x_i$ can appear in both normal or negated form; for example if the clause is $(x_i \vee \overline{x_j} \vee x_k)$ the corresponding NAE-SAT clause will be $(x_i \vee \overline{x_j} \vee x_k). (\overline{x_i} \vee x_j \vee \overline{x_k})$.

| NAE-SAT Clause $x = (x_i, x_j, x_k)$ | $E(\alpha_i, \alpha_j, \alpha_k)$ for a single clause ($\propto [T_{ij} + T_{jk} + T_{ki}]$) | | | |
|---|---|---|---|---|
| | $\alpha = (0,0,0) \equiv x = (1,1,1)$ | $\alpha = (0,0,\pi) \equiv x = (1,1,0)$ | $\alpha = (0,\pi,\pi) \equiv x = (1,0,0)$ | $\alpha = (\pi,\pi,\pi) \equiv x = (0,0,0)$ |
| $C_{NAE} = (x_i \vee x_j \vee x_k) . (\overline{x_i} \vee \overline{x_j} \vee \overline{x_k})$ | $\left(\frac{189}{256}\right)\pi A - \frac{3}{2}\pi A_s$ $C_{NAE} = 0$ | $\left(\frac{-51}{256}\right)\pi A - \frac{3}{2}\pi A_s$ $C_{NAE} = 1$ | $\left(\frac{-51}{256}\right)\pi A - \frac{3}{2}\pi A_s$ $C_{NAE} = 1$ | $\left(\frac{189}{256}\right)\pi A - \frac{3}{2}\pi A_s$ $C_{NAE} = 0$ |
| $C_{NAE} = (x_i \vee x_j \vee \overline{x_k}) . (\overline{x_i} \vee \overline{x_j} \vee x_k)$ | $\left(\frac{-51}{256}\right)\pi A - \frac{3}{2}\pi A_s$ $C_{NAE} = 1$ | $\left(\frac{189}{256}\right)\pi A - \frac{3}{2}\pi A_s$ $C_{NAE} = 0$ | $\left(\frac{-51}{256}\right)\pi A - \frac{3}{2}\pi A_s$ $C_{NAE} = 1$ | $\left(\frac{-51}{256}\right)\pi A - \frac{3}{2}\pi A_s$ $C_{NAE} = 1$ |
| $C_{NAE} = (x_i \vee \overline{x_j} \vee \overline{x_k}) . (\overline{x_i} \vee x_j \vee x_k)$ | $\left(\frac{-51}{256}\right)\pi A - \frac{3}{2}\pi A_s$ $C_{NAE} = 1$ | $\left(\frac{-51}{256}\right)\pi A - \frac{3}{2}\pi A_s$ $C_{NAE} = 1$ | $\left(\frac{189}{256}\right)\pi A - \frac{3}{2}\pi A_s$ $C_{NAE} = 0$ | $\left(\frac{-51}{256}\right)\pi A - \frac{3}{2}\pi A_s$ $C_{NAE} = 1$ |
| $C_{NAE} = (\overline{x_i} \vee \overline{x_j} \vee \overline{x_k}) . (x_i \vee x_j \vee x_k)$ | $\left(\frac{189}{256}\right)\pi A - \frac{3}{2}\pi A_s$ $C_{NAE} = 0$ | $\left(\frac{-51}{256}\right)\pi A - \frac{3}{2}\pi A_s$ $C_{NAE} = 1$ | $\left(\frac{-51}{256}\right)\pi A - \frac{3}{2}\pi A_s$ $C_{NAE} = 1$ | $\left(\frac{189}{256}\right)\pi A - \frac{3}{2}\pi A_s$ $C_{NAE} = 0$ |

**Fig. 2.** $E(\alpha_i, \alpha_j, \alpha_k)$ for a single NAE-3-SAT clause computed for different combination of the literals. It can be observed that the energy is minimum only when the NAE-SAT clause is satisfied. Only selected combinations have been shown here; a detailed table considering all combinations has been shown in Appendix IV.

To show that the energy corresponding to a clause, $[T_{ij} + T_{jk} + T_{ki}]$, is minimized when an NAE-3-SAT clause is satisfied, we consider the table in Fig. 2. It can be observed from the table that an NAE-SAT clause is satisfied only when $[T_{ij} + T_{jk} + T_{ki}]$ assumes the minimum value. Considering the inherent symmetry in the expression, only selected cases have been presented here. However, the complete table has been shown in Appendix IV.

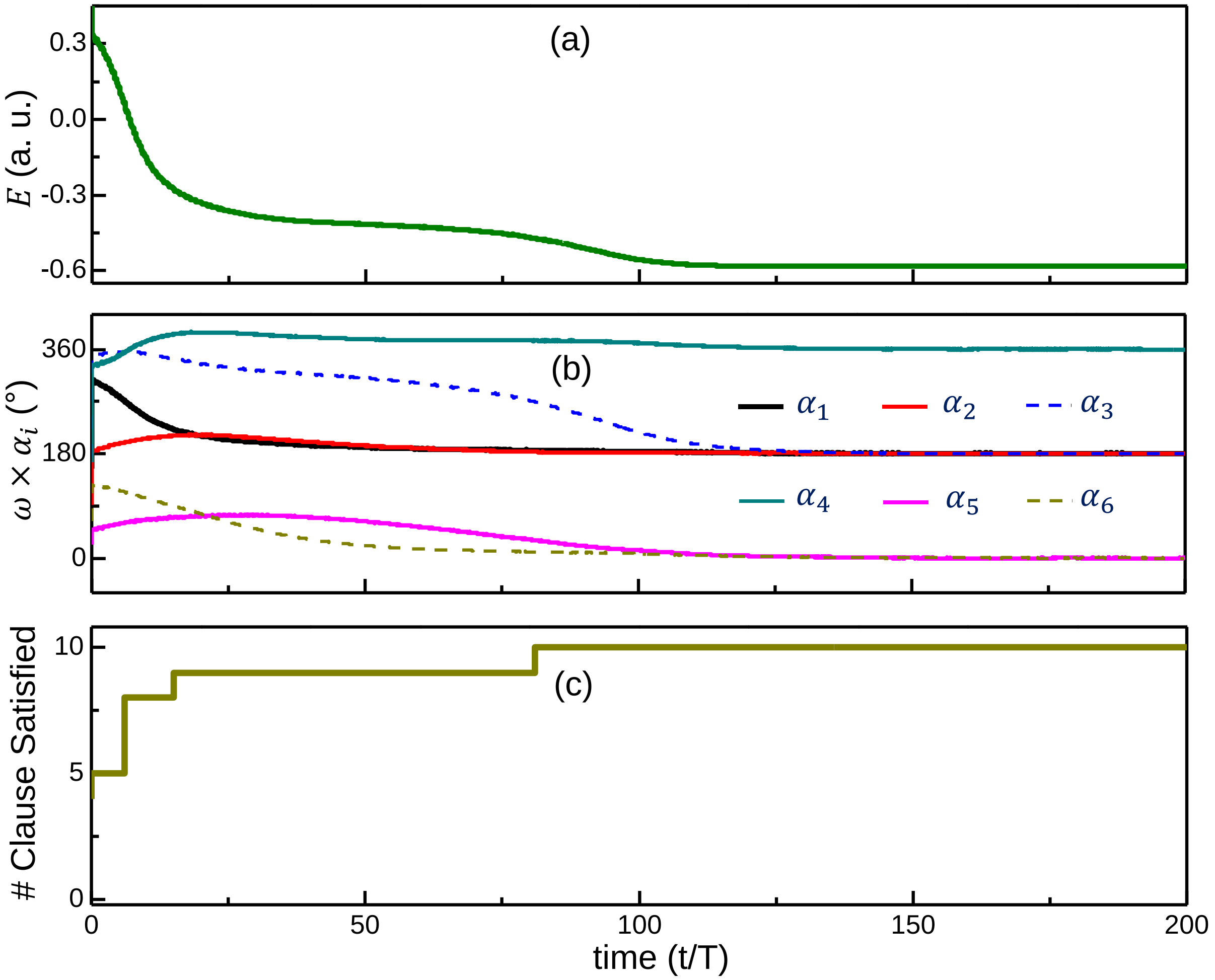

**Fig. 3.** Evolution of (a) $E$, (b) $\omega.\alpha_i$, and (c) number of satisfied clauses with time, for an illustrative problem with 6 variables and 10 clauses that is solved using the System II dynamics. The dynamics are obtained by simulating equation (10). In this simulation, $\omega = 2\pi$ is used such that T=1. Details of the SDE simulation are described in Appendix I.

Consequently, as the system evolves towards the global minimum of $E = \sum_{m=1}^{M} \beta_m\left(\alpha_i, \alpha_j, \alpha_k\right) + C - C_s$, it aims to maximize the number of satisfied NAE-3-SAT clauses (defined by $\left(x_i \vee x_j \vee x_k\right).\left(\overline{x_i} \vee \overline{x_j} \vee \overline{x_k}\right)$). In other words, it computes the solution to the Max-NAE-3-SAT problem. Figure 3 shows the solution for an illustrative NAE-3-SAT problem having 6 variables and 10 clauses. The oscillator dynamics are simulated using equation (10). However, equation (19) can also be used to compute the solution, as shown in Appendix V.

## Conclusion

In summary, we proposed two oscillator-inspired dynamical systems and demonstrated their ability to compute two different categories of the Boolean SAT problem. Our work helps advance the understanding of how oscillator-based systems can be designed and evaluated, specifically, for solving combinatorial optimization problems that have objective functions with degree greater than two. It also provides insights into how dynamical system formulation impacts the computing characteristics of the system. Further, our work reveals that, unlike problems with quadratic degree objective functions (e.g., MaxCut), oscillator models cannot be *directly* formulated as the difference of two oscillator phases. Consequently, the results presented here help expand the application of oscillator-based dynamical systems and the physics-inspired computing approach to a larger class of combinatorial optimization problems.

## Appendix I

Here, we describe the simulation approach used in the example problems considered in Fig. 1 and Fig. 3. We use a stochastic differential equation (SDE) solver that allows noise to be considered; noise can helps escape local minima in the phase space. The SDE that is used to evaluate the time evolution of oscillators' phases is given by,

$$d\alpha_{it} = f(\alpha_{it}, t)\, dt + \ a_n(\alpha_{it}, t) dw_t \qquad \text{(A1.1)}$$

where $f(\alpha_{it}, t)$ is the right-hand side of equation (8b) for System I and the right-hand side of equation (10) for System II, respectively; $a_n(\alpha_{it}, t)$ is the amplitude of noise. A time and phase independent noise amplitude of $5 \times 10^{-4}$ is used here. $w_t$ is a Wiener process [27]. We use the Runge-Kutta method of order four to develop the differential equation solver [28] in MATLAB. Values of $A$ and $A_s$ used in the simulation are:

| | $\boldsymbol{A}$ | $\boldsymbol{A_s}$ |
|---|---|---|
| System I | $\frac{10}{2\pi}$ | N/A |
| System II | $\frac{5}{2\pi}$ | $\frac{0.01}{2\pi}$ |

The illustrative 3-SAT / NAE-3-SAT problem used in all the above examples is given by:

$$\begin{aligned} Y \ =\ & (x_1 \vee x_2 \vee x_4) \wedge (x_1 \vee x_4 \vee x_5) \wedge (x_2 \vee x_3 \vee x_6) \wedge (\overline{x_2} \vee \overline{x_5} \vee \overline{x_6}) \\ & \wedge (x_1 \vee \overline{x_2} \vee x_6) \wedge (x_1 \vee x_4 \vee \overline{x_6}) \wedge (x_2 \vee \overline{x_4} \vee x_6) \wedge (\overline{x_1} \vee \overline{x_3} \vee \overline{x_4}) \qquad \text{(A1.2)} \\ & \wedge (\overline{x_4} \vee x_5 \vee \overline{x_6}) \wedge (\overline{x_3} \vee \overline{x_4} \vee \overline{x_5}) \end{aligned}$$

**Appendix II**

Here, we describe in detail the steps involved in the derivation of equation (12) for System II. The dynamics described by equation (11) are shown here (again) in:

$$\frac{d\alpha_i}{dt} = -\frac{A}{2}\left[\sum_{m=1}^{M}\left(c_{mi}\sin(t+\alpha_i)\left(\prod_{j=1;j\neq i}^{N}\left(\frac{1-c_{mj}\cos(t+\alpha_j)}{2}\right)\right)^2\right)\right.$$
$$\left.-\sum_{m=1}^{M}\left(\frac{1}{2}c_{mi}^2\sin\big(2(t+\alpha_i)\big)\left(\prod_{j=1;j\neq i}^{N}\left(\frac{1-c_{mj}\,cos(t+\alpha_j)}{2}\right)\right)^2\right)\right] \quad \text{(A2.1)}$$
$$-\sin(2t+2\alpha_i).\,[A_s\cos(2t)]$$

The $\cos(.)$ terms in equation (A2.1) are expressed by the following equation,

$$\left(\prod_{j=1;j\neq i}^{N}\left(\frac{1-c_{mj}\cos(t+\alpha_j)}{2}\right)\right)^2 \quad \text{(A2.2)}$$
$$=\sum_{\mu_N=-2}^{2}\cdots\sum_{\mu_2=-2}^{2}\sum_{\mu_1=-2}^{2}C_{\mu_1,\mu_2\ldots\mu_N;\neq\mu_i}\cos\left(\left(\sum_{j=1;\,j\neq i}^{N}|c_{mj}|\mu_j\right)t+\sum_{j=1;\,j\neq i}^{N}|c_{mj}|\mu_j\alpha_j\right)$$

Substituting (A2.2) in (A2.1),

$$\frac{d\alpha_i}{dt}$$

$$= -\frac{A}{2}\left[\sum_{m=1}^{M}\left(c_{mi}\sum_{\mu_N=-2}^{2}\cdots\sum_{\mu_2=-2}^{2}\sum_{\mu_1=-2}^{2}C_{\mu_1,\mu_2\ldots\mu_N;\neq\mu_i}\sin(t+\alpha_i)\,.\cos\left(\left(\sum_{j=1;\,j\neq i}^{N}|c_{mj}|\mu_j\right)t\right.\right.\right.$$

$$\left.\left.+\sum_{j=1;\,j\neq i}^{N}|c_{mj}|\mu_j\alpha_j\right)\right) \tag{A2.3}$$

$$-\sum_{m=1}^{M}\left(\frac{1}{2}c_{mi}^2\sum_{\mu_N=-2}^{2}\cdots\sum_{\mu_2=-2}^{2}\sum_{\mu_1=-2}^{2}C_{\mu_1,\mu_2\ldots\mu_N;\neq\mu_i}\sin\big(2(t+\alpha_i)\big).\cos\left(\left(\sum_{j=1;\,j\neq i}^{N}|c_{mj}|\mu_j\right)t\right.\right.$$

$$\left.\left.\left.+\sum_{j=1;\,j\neq i}^{N}|c_{mj}|\mu_j\alpha_j\right)\right)\right]-\sin(2t+2\alpha_i).\,[A_s\cos(2t)]$$

Equation (A2.3) can be written as,

$$\frac{d\alpha_i}{dt} = -\frac{A}{2}\left[\sum_{m=1}^{M}\left(c_{mi}\sum_{\mu_N=-2}^{2}\cdots\sum_{\mu_2=-2}^{2}\sum_{\mu_1=-2}^{2}C_{\mu_1,\mu_2\ldots\mu_N;\neq\mu_i}\,.\,\chi^{(1)}(t,\alpha_i)\,.\,B(t,\alpha)\right)\right.$$

$$\left.-\sum_{m=1}^{M}\left(\frac{1}{2}c_{mi}^2\sum_{\mu_N=-2}^{2}\cdots\sum_{\mu_2=-2}^{2}\sum_{\mu_1=-2}^{2}C_{\mu_1,\mu_2\ldots\mu_N;\neq\mu_i}\,.\,\chi^{(2)}(t,\alpha_i)\,.\,B(t,\alpha)\right)\right] \tag{A2.4}$$

$$-\sin(2t+2\alpha_i).\,[A_s\cos(2t)]$$

where, $\chi^{(1)}(t,\alpha_i)$ and $\chi^{(2)}(t,\alpha_i)$ are the first and the second harmonics of the perturbation projection vector (PPV) of the oscillator, respectively; $B(t,\alpha)$ is a perturbation which will have components ranging from the first harmonic to the $\left(2\sum_{j=1,j\neq i}^{N}|c_{mj}|\right)^{th}$ harmonic. As mentioned in the main text, we assume that the phase evolution happens on a much slower time scale than the oscillation frequency. Wang et al. [26] showed that such an equation (equation A2.4 here) can be formulated as a Multi-time Partial Differential

Equation (MPDE) and can be approximated by averaging over the fast time. The resulting approximation is essentially a cross-correlation of the PPV and the perturbation. Since cross correlation of the $i^{th}$ harmonic of the PPV with $j^{th}$ harmonic of the perturbation will be zero in cases where $i \neq j$, equation (A2.4) reduces to,

$$\frac{d\alpha_i}{dt} = -\frac{A}{2}\left[\sum_{m=1}^{M}\left(c_{mi}\sum_{\mu_N=-2}^{2}\cdots\sum_{\mu_2=-2}^{2}\sum_{\mu_1=-2}^{2}C_{\mu_1,\mu_2\ldots\mu_N;\neq\mu_i}\,.\,\chi^{(1)}(t,\alpha_i)\,.\,B^{(1)}(t,\alpha)\right)\right.$$
$$\left.-\sum_{m=1}^{M}\left(\frac{1}{2}c_{mi}^{2}\sum_{\mu_N=-2}^{2}\cdots\sum_{\mu_2=-2}^{2}\sum_{\mu_1=-2}^{2}C_{\mu_1,\mu_2\ldots\mu_N;\neq\mu_i}\,.\,\chi^{(2)}(t,\alpha_i)\,.\,B^{(2)}(t,\alpha)\right)\right] \tag{A2.5}$$
$$-\sin(2t+2\alpha_i).\,[A_s\cos(2t)]$$

Now, replacing $\chi$ and $B$ in (A2.5) with their original expressions from equation (A2.3), and by taking average with respect to the fast time, equation (A2.5) can be rewritten as,

$$\frac{d\alpha_i}{dt} = -A\sum_{m=1}^{M}\left(\sum_{\mu_N=-2}^{2}\ldots\sum_{\mu_2=-2}^{2}\sum_{\mu_1=-2}^{2}c_{mi}\;Q_1C^{1}_{\mu_1,\mu_2\ldots\mu_N;\neq\mu_i}\sin\left(\left(\alpha_i\right.\right.\right.$$
$$\left.\left.\left.-\sum_{j=1;\,j\neq i}^{N}\left|c_{mj}\right|\mu_j\alpha_j(t)\right|_{Q1}\right)\right)\right) \tag{A2.6}$$
$$+$$
$$A\sum_{m=1}^{M}\left(\sum_{\mu_N=-2}^{2}\ldots\sum_{\mu_2=-2}^{2}\sum_{\mu_1=-2}^{2}c_{mi}^{2}\;Q_2C^{2}_{\mu_1,\mu_2\ldots\mu_N;\neq\mu_i}\sin\left(2\alpha_i-\sum_{j=1;\,j\neq i}^{N}\left|c_{mj}\right|\mu_j\alpha_j(t)\right|_{Q2}\right)\right)$$
$$-A_{s1}\sin(2\alpha_i)$$

where, $Q_1 = 1$ when $\sum_{j=1;\ j\neq i}^{N} |c_{mj}| \mu_j = 1$ else $Q_1 = 0$; $Q_2 = 1$ when $\sum_{j=1;\ j\neq i}^{N} |c_{mj}| \mu_j = 2$ else $Q_2 = 0$. Equation (A2.6) is result shown in equation (12) in the main text.

**Appendix III**

Here, the energy expression (equation 18) for a clause containing 3 literals (System II) is derived. We assume that the 3 literals correspond to three distinct variables ($x_i$, $x_j$, $x_k$). Thus, equation (11) can be written as,

$$\frac{d\alpha_i}{dt} = -\frac{A}{2} 2^{-2N+6} \left[ \sum_{\substack{m=1; i\neq j\neq k; c_{mi}\neq 0 \\ c_{mj}\neq 0, c_{mk}\neq 0}}^{M} \left( c_{mi} \sin(t + \alpha_i) \left(\frac{1 - c_{mj}\cos(t+\alpha_j)}{2}\right)^2 \left(\frac{1 - c_{mk}\cos(t+\alpha_k)}{2}\right)^2 \right) - \sum_{\substack{m=1; i\neq j\neq k; c_{mi}\neq 0 \\ c_{mj}\neq 0, c_{mk}\neq 0}}^{M} \left( \frac{1}{2} c_{mi}^2 \sin(2t + 2\alpha_i) \left(\frac{1 - c_{mj}\cos(t+\alpha_j)}{2}\right)^2 \left(\frac{1 - c_{mk}\cos(t+\alpha_k)}{2}\right)^2 \right) \right] - \sin(2t + 2\alpha_i) \left[A_s \cos(2t)\right] \tag{A3.1}$$

Equation (A3.1) can be simplified to,

$$\frac{d\alpha_i}{dt} = -\frac{A}{2}2^{-2N+2}\Bigg[\sum_{\substack{m=1;i\neq j\neq k;c_{mi}\neq 0\\ c_{mj}\neq 0,c_{mk}\neq 0}}^{M}\Bigg(c_{mi}\sin(t+\alpha_i)\left(1+\frac{1}{2}c_{mj}^2 - 2c_{mj}\cos(t+\alpha_j)\right.$$

$$\left.+\frac{1}{2}c_{mj}^2\cos(2t+2\alpha_j)\right)\left(1+\frac{1}{2}c_{mk}^2 - 2c_{mk}\cos(t+\alpha_k)\right.$$

$$\left.\left.+\frac{1}{2}c_{mk}^2\cos(2t+2\alpha_k)\right)\right) \tag{A3.2}$$

$$-\sum_{\substack{m=1;i\neq j\neq k;c_{mi}\neq 0\\ c_{mj}\neq 0,c_{mk}\neq 0}}^{M}\left(\frac{1}{2}c_{mi}^2\sin(2t+2\alpha_i)\left(1+\frac{1}{2}c_{mj}^2 - 2c_{mj}\cos(t+\alpha_j)\right.\right.$$

$$\left.+\frac{1}{2}c_{mj}^2\cos(2t+2\alpha_j)\right)\left(1+\frac{1}{2}c_{mk}^2 - 2c_{mk}\cos(t+\alpha_k)\right.$$

$$\left.\left.+\frac{1}{2}c_{mk}^2\cos(2t+2\alpha_k)\right)\right)\Bigg] - \sin(2t+2\alpha_i).\,[A_s\cos(2t)]$$

After subsequent simplification of the terms and averaging over the fast time, equation (A3.2) can be expressed as,

$$\frac{d\alpha_i}{dt} = -\frac{\pi A}{2}2^{-2N+2}\sum_{\substack{m=1;i\neq j\neq k;c_{mi}\neq 0\\ c_{mj}\neq 0,c_{mk}\neq 0}}^{M}\left[-2c_{mi}c_{mj}\left(1+\frac{1}{2}c_{mk}^2\right)\sin(\alpha_i-\alpha_j)\right.$$

$$-2c_{mi}c_{mk}\left(1+\frac{1}{2}c_j^2\right)\sin(\alpha_i-\alpha_k) - \frac{1}{2}c_{mi}c_{mj}c_{mk}^2\sin(\alpha_i+\alpha_j-2\alpha_k)$$

$$\left.-\frac{1}{2}c_{mi}c_{mk}c_{mj}^2\sin(\alpha_i+\alpha_k-2\alpha_j)\right]$$

$$+\frac{\pi A}{2}2^{-2N+2}\sum_{\substack{m=1;i\neq j\neq k;c_{mi}\neq 0\\ c_{mj}\neq 0,c_{mk}\neq 0}}^{M}\left[\frac{1}{4}c_{mi}^2c_{mk}^2\left(1+\frac{1}{2}c_{mj}^2\right)\sin(2\alpha_i-2\alpha_k)\right. \tag{A3.3}$$

$$+\frac{1}{4}c_{mi}^2c_{mj}^2\left(1+\frac{1}{2}c_{mk}^2\right)\sin(2\alpha_i-2\alpha_j)$$

$$\left.+c_{mi}^2c_{mj}c_{mk}\sin(2\alpha_i-\alpha_j-\alpha_k)\right] - \pi A_s\sin(2\alpha_i)$$

Equation (A3.3) can be further simplified as

$$\frac{d\alpha_i}{dt} = \frac{\pi A}{2} 2^{-2N+2} \sum_{\substack{m=1; i\neq j\neq k; c_{mi}\neq 0 \\ c_{mj}\neq 0, c_{mk}\neq 0}}^{M} \Big[ 2c_{mi}c_{mj}\left(1+\frac{1}{2}c_{mk}^2\right)\sin(\alpha_i-\alpha_j)$$
$$+ 2c_{mi}c_{mk}\left(1+\frac{1}{2}c_{mj}^2\right)\sin(\alpha_i-\alpha_k) + \frac{1}{2}c_{mi}c_{mj}c_{mk}^2\sin(\alpha_i+\alpha_j-2\alpha_k)$$
$$+ \frac{1}{2}c_{mi}c_{mk}c_{mj}^2\sin(\alpha_i+\alpha_k-2\alpha_j) + \frac{1}{4}c_{mi}^2c_{mk}^2\left(1+\frac{1}{2}c_{mj}^2\right)\sin(2\alpha_i-2\alpha_k) \qquad \text{(A3.4)}$$
$$+ \frac{1}{4}c_{mi}^2c_{mj}^2\left(1+\frac{1}{2}c_{mk}^2\right)\sin(2\alpha_i-2\alpha_j)$$
$$+ c_{mi}^2c_{mj}c_{mk}\sin(2\alpha_i-\alpha_j-\alpha_k)\Big] - \pi A_s\sin(2\alpha_i)$$

Using equation (15), the Lyapunov function for the above dynamics can be defined as,

$$E(\alpha) = \pi A.2^{-2N+1} \sum_{i=1}^{N} \sum_{\substack{m=1; i\neq j\neq k; c_{mi}\neq 0 \\ c_{mj}\neq 0, c_{mk}\neq 0}}^{M} \Big[ 2c_{mi}c_{mj}\left(1+\frac{1}{2}c_{mk}^2\right)\cos(\alpha_i-\alpha_j)$$
$$+ 2c_{mi}c_{mk}\left(1+\frac{1}{2}c_{mj}^2\right)\cos(\alpha_i-\alpha_k) + \frac{1}{2}c_{mi}c_{mj}c_{mk}^2\cos(\alpha_i+\alpha_j-2\alpha_k)$$
$$+ \frac{1}{2}c_{mi}c_{mk}c_{mj}^2\cos(\alpha_i+\alpha_k-2\alpha_j) + \frac{1}{8}c_{mi}^2c_{mk}^2\left(1+\frac{1}{2}c_{mj}^2\right)\cos(2\alpha_i-2\alpha_k) \qquad \text{(A3.5)}$$
$$+ \frac{1}{8}c_{mi}^2c_{mj}^2\left(1+\frac{1}{2}c_{mk}^2\right)\cos(2\alpha_i-2\alpha_j) + \frac{1}{2}c_{mi}^2c_{mj}c_{mk}\cos(2\alpha_i-\alpha_j-\alpha_k)\Big]$$
$$- \sum_{i=1}^{N}\frac{\pi A_s}{2}\cos(2\alpha_i)$$

Equation (A3.5) expresses the energy function for the NAE-3-SAT problem. Equation (A3.5) is same as the result shown in equation (18)

## Appendix IV

Here, we show that an NAE-3-SAT clause is satisfied only when the corresponding energy term associated with the clause is minimized (all cases).

| NAE-SAT Clause $x = (x_i, x_j, x_k)$ | $E(\alpha_i, \alpha_j, \alpha_k)$ for a single clause ($\propto [T_{ij} + T_{jk} + T_{ki}]$) | | | | | | | |
|---|---|---|---|---|---|---|---|---|
| | $(0,0,0) \equiv x = (1,1,1)$ | $(0,0,\pi) \equiv x = (1,1,0)$ | $(0,\pi,0) \equiv x = (1,0,1)$ | $(\pi,0,0) \equiv x = (0,1,1)$ | $(\pi,\pi,0) \equiv x = (0,0,1)$ | $(0,\pi,\pi) \equiv x = (1,0,0)$ | $(\pi,0,\pi) \equiv x = (0,1,0)$ | $(\pi,\pi,\pi) \equiv x = (0,0,0)$ |
| $C_{NAE} = (x_i \vee x_j \vee x_k) . (\overline{x_i} \vee \overline{x_j} \vee \overline{x_k})$ | $E_1$ $C_{NAE} = 0$ | $E_2$ $C_{NAE} = 1$ | $E_2$ $C_{NAE} = 1$ | $E_2$ $C_{NAE} = 1$ | $E_2$ $C_{NAE} = 1$ | $E_2$ $C_{NAE} = 1$ | $E_2$ $C_{NAE} = 1$ | $E_1$ $C_{NAE} = 0$ |
| $C_{NAE} = (x_i \vee x_j \vee \overline{x_k}) . (\overline{x_i} \vee \overline{x_j} \vee x_k)$ | $E_2$ $C_{NAE} = 1$ | $E_1$ $C_{NAE} = 0$ | $E_2$ $C_{NAE} = 1$ | $E_2$ $C_{NAE} = 1$ | $E_1$ $C_{NAE} = 0$ | $E_2$ $C_{NAE} = 1$ | $E_2$ $C_{NAE} = 1$ | $E_2$ $C_{NAE} = 1$ |
| $C_{NAE} = (x_i \vee \overline{x_j} \vee x_k) . (\overline{x_i} \vee x_j \vee \overline{x_k})$ | $E_2$ $C_{NAE} = 1$ | $E_2$ $C_{NAE} = 1$ | $E_1$ $C_{NAE} = 0$ | $E_2$ $C_{NAE} = 1$ | $E_2$ $C_{NAE} = 1$ | $E_2$ $C_{NAE} = 1$ | $E_1$ $C_{NAE} = 0$ | $E_2$ $C_{NAE} = 1$ |
| $C_{NAE} = (\overline{x_i} \vee x_j \vee x_k) . (x_i \vee \overline{x_j} \vee \overline{x_k})$ | $E_2$ $C_{NAE} = 1$ | $E_2$ $C_{NAE} = 1$ | $E_2$ $C_{NAE} = 1$ | $E_1$ $C_{NAE} = 0$ | $E_2$ $C_{NAE} = 1$ | $E_1$ $C_{NAE} = 0$ | $E_2$ $C_{NAE} = 1$ | $E_2$ $C_{NAE} = 1$ |
| $C_{NAE} = (\overline{x_i} \vee \overline{x_j} \vee x_k) . (x_i \vee x_j \vee \overline{x_k})$ | $E_2$ $C_{NAE} = 1$ | $E_1$ $C_{NAE} = 0$ | $E_2$ $C_{NAE} = 1$ | $E_2$ $C_{NAE} = 1$ | $E_1$ $C_{NAE} = 0$ | $E_2$ $C_{NAE} = 1$ | $E_2$ $C_{NAE} = 1$ | $E_2$ $C_{NAE} = 1$ |
| $C_{NAE} = (x_i \vee \overline{x_j} \vee \overline{x_k}) . (\overline{x_i} \vee x_j \vee x_k)$ | $E_2$ $C_{NAE} = 1$ | $E_2$ $C_{NAE} = 1$ | $E_2$ $C_{NAE} = 1$ | $E_1$ $C_{NAE} = 0$ | $E_2$ $C_{NAE} = 1$ | $E_1$ $C_{NAE} = 0$ | $E_2$ $C_{NAE} = 1$ | $E_2$ $C_{NAE} = 1$ |
| $C_{NAE} = (\overline{x_i} \vee x_j \vee \overline{x_k}) . (x_i \vee \overline{x_j} \vee x_k)$ | $E_2$ $C_{NAE} = 1$ | $E_2$ $C_{NAE} = 1$ | $E_1$ $C_{NAE} = 0$ | $E_2$ $C_{NAE} = 1$ | $E_2$ $C_{NAE} = 1$ | $E_2$ $C_{NAE} = 1$ | $E_1$ $C_{NAE} = 0$ | $E_2$ $C_{NAE} = 1$ |
| $C_{NAE} = (\overline{x_i} \vee \overline{x_j} \vee \overline{x_k}) . (x_i \vee x_j \vee x_k)$ | $E_1$ $C_{NAE} = 0$ | $E_2$ $C_{NAE} = 1$ | $E_2$ $C_{NAE} = 1$ | $E_2$ $C_{NAE} = 1$ | $E_2$ $C_{NAE} = 1$ | $E_2$ $C_{NAE} = 1$ | $E_2$ $C_{NAE} = 1$ | $E_1$ $C_{NAE} = 0$ |

$E_1 = \frac{189}{256}\pi A - \frac{3}{2}\pi A_s$ $\qquad$ $E_2 = \frac{-51}{256}\pi A - \frac{3}{2}\pi A_s$

**Fig. 4.** $E(\alpha_i, \alpha_j, \alpha_k)$ corresponding to a single NAE-3-SAT clause computed for all the possible combinations of the literals and phases. Here, $E_1 > E_2$.

It can be observed from the table in Fig. 4 that an NAE-3-SAT clause is satisfied only when the corresponding energy term associated with the clause is minimized. Consequently, the decreasing nature of the energy function ensures that the system evolves towards a state that maximizes the number of satisfied NAE-3-SAT clauses.

**Appendix V**

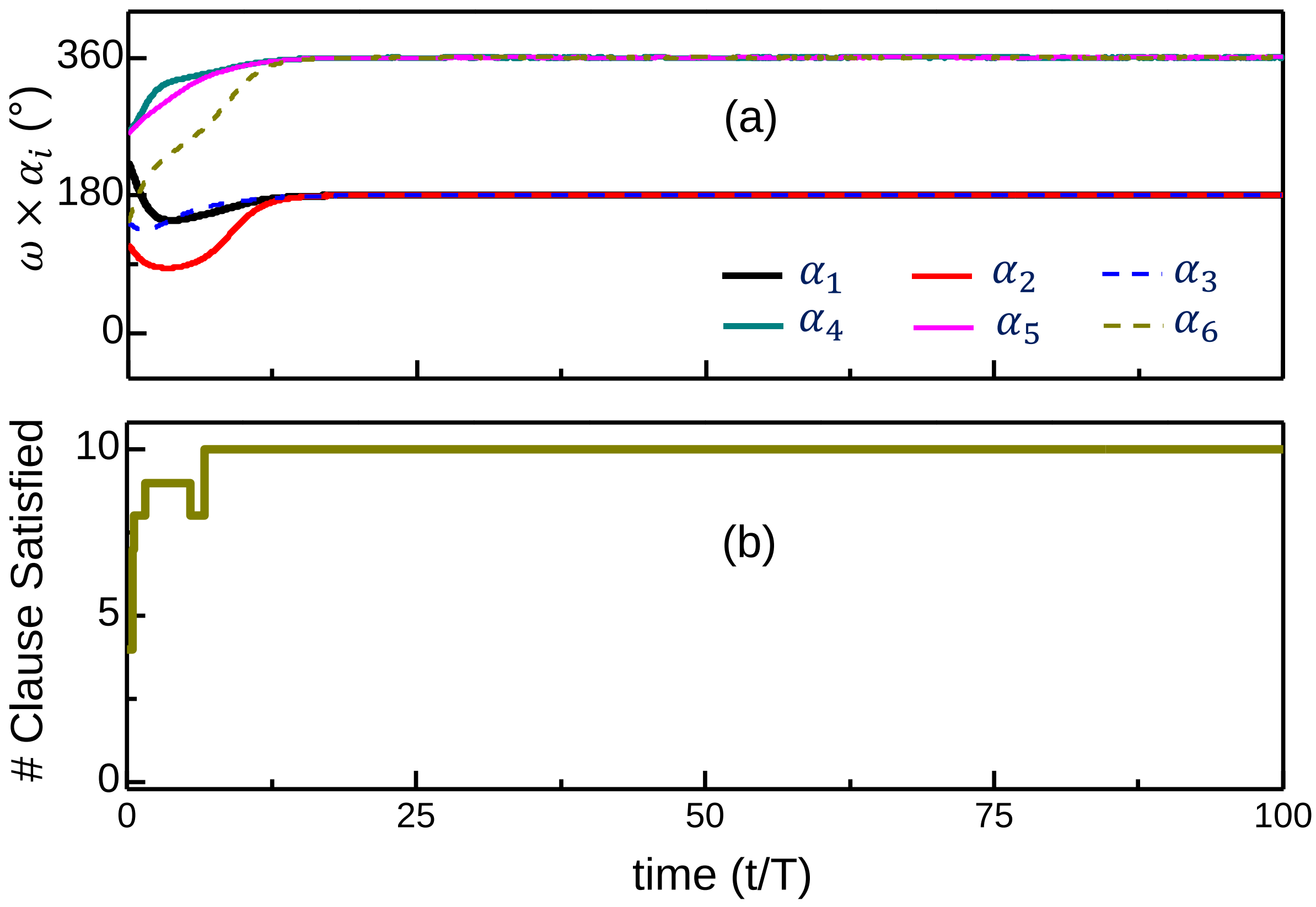

**Fig. 5.** Evolution of (a) oscillator phases, calculated using equation (19); and (b) No. of satisfied clauses, with time. $\omega = 2\pi$ is used such that $T = 1$. The problem considered here is the same as that considered in Fig. 3.

We show in Fig. 5 the NAE-3-SAT solution can also be computed using equation (19).

**Acknowledgment:**

This work was supported by NSF ASCENT grant (No. 2132918). We would like to thank Professor Avik Ghosh from University of Virginia for insightful discussions.